\documentclass[a4paper, 10pt, oneside]{article}
\usepackage{amsmath}
\usepackage{amssymb}
\usepackage{verbatim}
\usepackage{setspace}

\usepackage[amsmath,thmmarks]{ntheorem}

\theorembodyfont{\normalfont}
\newtheorem{Thm}[equation]{Theorem}

\newtheorem{Lem}[equation]{Lemma}

\newtheorem{Cor}[equation]{Corollary}

\theoremstyle{nonumberplain} \theoremheaderfont{\normalfont\itshape}
\theorembodyfont{\normalfont} \theoremseparator{.}
\theoremsymbol{\ensuremath{\square}}

\begin{document}
\begin{center}
\textbf{{\sc Symmetric Presentations of Coxeter Groups}}
\bigskip

\textbf{{\sc Ben Fairbairn}} - \texttt{bfairbairn@ems.bbk.ac.uk}

\bigskip

\emph{Department of Economics, Mathematics and Statistics, Birkbeck,\\ University of London, Malet Street, London WC1E 7HX}
\end{center}

\begin{quote}
\begin{center}
Abstract
\end{center}
We apply the techniques of symmetric generation to establish the
standard presentations of the finite simply laced irreducible finite
Coxeter groups, that is the Coxeter groups of types $A_n$, $D_n$ and
$E_n$, and show that these are naturally arrived at purely through
consideration of certain natural actions of symmetric groups. We go
on to use these techniques to provide explicit representations of
these groups.
\end{quote}


\section{Introduction}
A \emph{Coxeter diagram} of a presentation is a graph in which the
vertices correspond to involutory generators and an edge is labeled
with the order of the product of its two endpoints. Commuting
vertices are not joined and an edge is left unlabeled if the
corresponding product has order three. A Coxeter diagram and its
associated group are said to be \emph{simply laced} if all the edges
of the graph are unlabeled. In \cite{ATLASconference} Curtis notes
that if such a diagram has a ``tail'' of length at least two, as in
Figure I, then we see that the generator corresponding to the
terminal vertex, $a_r$, commutes with the subgroup generated by the
subgraph $\mathcal{G}_0$.
\bigskip

\setlength{\unitlength}{1mm}
\begin{picture}(0, 0)
\put(30,0){\circle{20}}\put(37,0){\line(1,0){40}}\multiput(57,0)(20,0){2}{\circle*{2}}
\put(28.5,-1){$\mathcal{G}_0$}\put(56,-3){$a_{r-1}$}\put(76,-3){$a_r$}
\end{picture}
\bigskip\bigskip
\begin{center}
Figure I: A Coxeter diagram with a tail.
\end{center}

\bigskip

 In this paper we slightly generalize the notion of a ``graph with a tail"
 and in doing so
 provide symmetric presentations for all the simply laced
irreducible finite Coxeter groups with the aid of little more than a
single short relation. These in turn readily give rise to natural
representations of these groups.

Presentations of groups having certain types of symmetry properties
have been considered at least since Coxeter's work \cite{CoxOrig} in
1959 and have proved useful not only in providing natural and
elementary definitions of groups but have also been of great
computational use. In \cite{CurtisFairbairn} Curtis and the author
used one kind of symmetric presentation for the Conway group
$\cdot$0 obtained by Bray and Curtis in \cite{revisited} to
represent elements of $\cdot$0 as a string of at most 64 symbols and
typically far fewer. This represents a considerable saving compared
to representing an element of $\cdot$0 as a permutation of 196560
symbols or as a $24 \times 24$ matrix (ie as a string of 24$^2$=576
symbols). More in depth discussions of symmetric generation more generally may be found in \cite{Curtismain,ATLASconference,Survey}.

The presentations given here, whilst not new, do provide an
excellent example of how the techniques of symmetric generation may
be used to arrive at very natural constructions of groups, and in
seeing how these presentations may in turn lead to highly symmetric
representations of these groups. Whilst recent results of the author
and M\"{u}ller \cite{muller} generalize our Main Theorem to a wider
class of Coxeter groups, the symmetric presentations there are not
well motivated (indeed it is the results presented here that provide
the main motivation for the results of \cite{muller}); may not be
arrived at as naturally as those presented here are and do not easily lead
to explicit representations (the matrices we are naturally lead to for the representations of the groups considered here being strikingly simple in nature).

For the basic definitions and notation for Coxeter groups used
throughout this paper we refer the reader to the book of Humphreys \cite{Coxetergroups}. Throughout we shall use the standard {\sc Atlas} notation
for groups found in \cite{ATLAS}.

This article is organised as follows. In Section 2 we outline the basic techniques of involutory symmetric
generation. In Section 3 we state our main theorem and the barriers
to further extension. In Section \ref{mot} we will show how general
results in symmetric generation naturally lead us straight to the
presentations considered in this paper. In Section \ref{enumerate}
we perform a coset enumeration necessary to prove our main theorem.
In Section \ref{reps} we use the symmetric presentations of the main
theorem to construct real representations of the groups concerned
and in doing so complete the proof. In Section \ref{modrep} we
construct $\mathbb{Z}_2$-representations from our real
representations in the $E_n$ cases to identify these groups as
$\Bbb{Z}_2$ matrix groups.

\section{Involutory Symmetric Generation}

We shall describe here only the case when the symmetric generators
are involutions as originally discussed by Bray, Curtis and Hammas
in \cite{involutorygenerators}. For a discussion of the more general
case see \cite[Section III]{Curtismain}.

Let $2^{\star n}$ denote the free product of $n$ involutions. We
write $\lbrace t_1,\dots,t_n\rbrace$ for a set of generators of this
free product. A permutation $\pi\in S_n$ induces an automorphism of
this free product $\hat\pi$ by permuting its generators, ie
$t_i^{\hat\pi}=t_{\pi(i)}$. Given a subgroup $N\leq S_n$ we can form
a semi-direct product $\cal{P}$=2$^{\star n}\colon N$ where, for
$\pi\in N$, $\pi^{-1}t_i\pi=t_{\pi(i)}$. When $N$ is transitive we
call $\cal{P}$ a \emph{progenitor}. We call $N$ the \emph{control
group} of $\cal{P}$ and the $t_i$'s the \emph{symmetric generators}.
Elements of $\cal{P}$ can all be written in the form $\pi w$ with
$\pi\in N$ and $w$ is a word in the symmetric generators, so any
homomorphic image of the progenitor can be obtained by factoring out
relations of the form $\pi w=1$. We call such a homomorphic image
that is finite a \emph{target group}. If $G$ is the target group
obtained by factoring the progenitor 2$^{\star n}\colon N$ by the
relators $\pi_1w_1,\pi_2w_2,\ldots$ we write
\[\frac{2^{\star n}\colon N}{\pi_1w_1,\pi_2w_2,\ldots}\cong G.\]
In keeping with the now traditional notational conventions used in works discussing symmetric generation, 
we write $N$ both for the control group and its image in $G$ and
refer to both simply as `the control group'. Similarly we shall
write $t_i$ both for a symmetric generator and its image in $G$ and
we shall refer to both as a `symmetric generator'.

To decide whether a given homomorphic image of a progenitor is
finite, we shall perform a coset enumeration. Given a word in the
symmetric generators, $w$, we define the {\em coset stabilizing
subgroup} of the coset $Nw$ to be the subgroup
\[ N^{(w)}:=\lbrace\pi\in N | Nw\pi=Nw\rbrace\leq N. \]
This is clearly a subgroup of $N$ and there are $|N:N^{(w)}|$ right
cosets of $N^{(w)}$ in $N$ contained in the double coset $NwN\subset
G$. We will write [$w$] for the double coset $NwN$ and [$\star$]
will denote the coset $[id_N]=N$. We shall write $w\sim w^\prime$ to
mean $[w]=[w^\prime]$. We can enumerate these cosets using
procedures such as the Todd-Coxeter algorithm, which can readily be
programmed into a computer. The sum of the numbers $|N:N^{(w)}|$
then gives the index of $N$ in $G$, and we are thus able to
determine the order of $G$ and in doing so prove it is finite.

In particular if the target group corresponds to the group defined
by a Coxeter diagram with a tail, then removing the vertex at the
end of the tail provides a control group for a symmetric
presentation with the vertex itself acting as a symmetric generator.

A family of results suggest that this approach lends itself to the
construction of groups with low index perfect subgroups. For
instance:

\begin{Lem}If $N$ is perfect and primitive, then
$|\cal{P}\colon\cal{P}'|$=2 and $\cal{P}''=\cal{P}'$.\end{Lem}
\begin{Cor}\label{perfect} If $N$ is perfect and primitive then any image of $\cal{P}$ possesses a perfect subgroup of
index at most 2. In particular any homomorphic image of $\cal{P}$
satisfying a relation of odd length is perfect.\end{Cor}

For proofs of these results see \cite[Theorem 1, p.356]{CurtisJ1}.

The next lemma, whilst easy to state and prove, has turned out to be
extremely powerful in leading to constructions of groups in terms of
symmetric generating sets, most notably a majority of the sporadic
simple groups \cite{Curtismain}.

\begin{Lem}\label{famous}
$$\langle t_i, t_j \rangle\cap N\leq C_N(Stab_N(i,j)).$$
\end{Lem}

Given a pair of symmetric generators $t_1$ and $t_2$, Lemma
\ref{famous} tells us which permutations $\pi\in N$ may be written
as a word in $t_1$ and $t_2$ but gives us no indication of the
length of such a word. Naturally we wish to factor a given
progenitor by the shortest and most easily understood relation
possible. The following lemma shows that in many circumstances a
relation of the form $\pi t_1t_2t_1$ is of great interest.

\begin{Lem}\label{primlem}
Let $G=\langle\mathcal{T}\rangle$, where
$\mathcal{T}=\{t_1,\ldots,t_n\}\subseteq G$ is a set of involutions
in $G$ with $N=N_G(\mathcal{T})$ acting primitively on $\mathcal{T}$
by conjugation. (Thus $G$ is a homomorphic image of the progenitor
$2^{\star n}:N$.) If $t_1t_i\in N$, $t_1\not\in N$ for some
$i\not=1$, then $|G|=2|N|$.
\end{Lem}

For proofs of these results see \cite[p.58 and p.59]{Curtismain}.

\section{The Main Theorem}

In the notation of the last section we will prove:

\begin{Thm} \label{main}Let $S_n$ be the symmetric group acting on $n$
objects and $W(\Phi)$ denote the Weyl group of the root system
$\Phi$. Then:

\begin{enumerate}
\item
\[ \frac{2^{\star {n \choose 1}}:S_n}{(t_1(12))^3}\cong W(\mbox{A}_{n}) \]
\item
\[ \frac{2^{\star {n \choose 2}}:S_n}{(t_{12}(23))^3}\cong W(\mbox{D}_n) \mbox{ for } n\geq4\]
\item
\[ \frac{2^{\star {n \choose 3}}:S_n}{(t_{123}(34))^3}\cong W(\mbox{E}_{n}) \mbox{ for }
n=6,7,8.\]
\end{enumerate}
\end{Thm}

In case (1) the action of $S_n$ defining the progenitor is the
natural action of $S_n$ on $X:=\{1,\ldots,n\}$; in case (2) the
action of $S_n$ defining the progenitor is the action of $S_n$ on
the 2-element subsets of $X$ and in case (3) the action of $S_n$
defining the progenitor is the action of $S_n$ on the 3-element
subsets of $X$.

Case (1) of Theorem \ref{main} has been noted by various authors
before \cite[Theorem 3.2, p.63]{Curtismain}, but we include it here
for completeness.

More suggestively we can express these symmetric presentations as
Coxeter diagrams as given in Figure II. (Notice that from the
presentations given in this Theorem, without even drawing any
Coxeter diagrams, the exceptional coincidences of D$_3$=A$_3$ and
E$_5$=D$_5$ are immediate since ${3 \choose 2}={3 \choose 1}$ and
${5 \choose 3}={5\choose 2}$).

We remark that the natural pattern of applying the relation
$(t_{1,\ldots,k}(k,k+1))^3$ to the progenitor $2^{\star {n \choose
k}}:S_n$ to produce a finite image does not extend further. In
\cite{Fischer}, Bray, Curtis, Parker and Wiedorn, prove the
symmetric presentation:
\[\frac{2^{\star {n \choose 4}}:\mbox{S}_8}{(t_{1234}(45))^3,t_{1234}t_{5678}}\cong
W(\mbox{E}_7)\cong \mbox{S}_6(2)\times2.\] The second relation,
which simply identifies a 4-element subset with its complement so
that the symmetric generators correspond to partitions of the eight
points into two fours, is necessary for the coset enumeration to
terminate; hence the pattern does not continue when the control
group is the full symmetric group. However, using a control group
smaller than the full symmetric group can resolve this problem. In
\cite{revisited} Bray and Curtis prove that:
\[\frac{2^{\star{24\choose4}}\colon\mbox{M}_{24}}{\pi t_{ab}t_{ac}t_{ad}}\cong\cdot0,\]
where M$_{24}$ denotes the largest of the sporadic simple Mathieu
groups; $a,b,c$ and $d$ are pairs of points the union of which is a
block of the $\mathcal{S}$(5,8,24) Steiner system on which M$_{24}$
naturally acts (see the {\sc Atlas}, \cite[p.94]{ATLAS}); $\cdot0$
is the full cover group of the largest sporadic simple Conway group
(see the {\sc Atlas}, \cite[p.180]{ATLAS}) and $\pi\in$ M$_{24}$ is
the unique permutation of M$_{24}$ set-wise fixing the sextets
defined by each of the symmetric generators whose use is motivated
by Lemma \ref{famous}.

\begin{figure}
\bigskip
\hspace{45mm}{\makebox[60pt]{%
\hspace{30mm}A$_n$:\hspace{65mm} \unitlength=1pt
\begin{picture}(20,20)(-30,-10)
\put(-200,5){\line(1,0){70}} \multiput(-200,5)(30,0){3}{\circle*{4}}
\put(-200,5){\line(0,-1){30}} \put(-200,-25){\circle*{4}}
\multiput(-129,5)(5,0){8}{\line(1,0){2}}
\put(-100,5){\line(1,0){40}} \multiput(-90,5)(30,0){2}{\circle*{4}}
\put(-210,12){(1,2)\hspace{3mm}(2,3)\hspace{3mm}(3,4)\hspace{19mm}($n-1$,$n$)}
\put(-195,-28){$t_1$}
\end{picture}%
}%
}

\bigskip\bigskip\bigskip
\hspace{45mm}{\makebox[60pt]{%
\hspace{30mm}D$_n$:\hspace{65mm} \unitlength=1pt
\begin{picture}(20,20)(-30,-10)
\put(-200,5){\line(1,0){70}} \multiput(-200,5)(30,0){3}{\circle*{4}}
\put(-170,5){\line(0,-1){30}} \put(-170,-25){\circle*{4}}
\multiput(-129,5)(5,0){8}{\line(1,0){2}}
\put(-100,5){\line(1,0){40}} \multiput(-90,5)(30,0){2}{\circle*{4}}
\put(-210,12){(1,2)\hspace{3mm}(2,3)\hspace{3mm}(3,4)\hspace{19mm}($n-1$,$n$)}
\put(-165,-28){$t_{1,2}$}
\end{picture}%
}%
}

\bigskip\bigskip\bigskip
\hspace{45mm}{\makebox[60pt]{%
\hspace{30mm}E$_n$:\hspace{65mm}\unitlength=1pt
\begin{picture}(20,20)(-30,-10)
\put(-200,5){\line(1,0){70}} \multiput(-200,5)(30,0){3}{\circle*{4}}
\put(-140,5){\line(0,-1){30}} \put(-140,-25){\circle*{4}}
\multiput(-129,5)(5,0){8}{\line(1,0){2}}
\put(-100,5){\line(1,0){40}} \multiput(-90,5)(30,0){2}{\circle*{4}}
\put(-210,12){(1,2)\hspace{3mm}(2,3)\hspace{3mm}(3,4)\hspace{19mm}($n-1$,$n$)}
\put(-135,-28){$t_{1,2,3}$}
\end{picture}%
}%
}\bigskip\bigskip\label{CoxDiags}

\begin{center}
Figure II: Symmetric presentations as Coxeter diagrams.
\end{center}
\end{figure}

The proof of Theorem \ref{main} is as follows. In Section
\ref{enumerate} we enumerate the double cosets $NwN$ in each case to
verify that the orders of the target groups are at most the orders
claimed in Theorem \ref{main}. In Section \ref{reps} we exhibit
elements of the target groups that generate them and satisfy the
additional relations, thereby providing lower bounds for the orders
and verifying the presentations.

\section{Motivating the Relations of Theorem \ref{main}}\label{mot}

In this section we will show how the relators used in Theorem
\ref{main} may be arrived at naturally by considering the natural
actions of the control group used to define the progenitors
appearing in the Main Theorem.

Given Lemma \ref{famous} it is natural to want to
compute $C_{S_n}(Stab_{S_n}(1,2))$. In the $A_n$ case we find
\[ Stab_{S_n}(1,2) =\left\{ \begin{array}{ll}
         \langle id\rangle & \mbox{if }n\in\{2,3\};\\
         \langle(3,4),(3,\ldots,n)\rangle & \mbox{if }n\geq4.\end{array} \right. \]

\noindent calculating $C_{S_n}(Stab_{S_n}(1,2))$ thus gives us

\[ C_{S_n}(Stab_{S_n}(1,2)) =\left\{ \begin{array}{ll}
\langle (1,2)\rangle & \mbox{if }n=2\mbox{ or }n\geq5;\\
         \langle (1,2),(1,2,3)\rangle & \mbox{if }n=3;\\
         \langle (1,2),(3,4)\rangle & \mbox{if }n=4.\end{array} \right. \]

For $n\geq5$ we see that $\langle t_1,t_2\rangle\cap N\leq\langle
(1,2)\rangle$. Lemma \ref{primlem} now tells us that the shortest
natural relator worth considering is $(1,2)t_1t_2t_1$ which we
rewrite more succinctly as $(t_1(12))^3$. We are thus naturally led
to considering the factored progenitor

\[ \frac{2^{\star {n \choose 1}}:S_n}{(t_1(12))^3}.\]

 Recall that
$S_n$ is the symmetric group acting on $n$ objects. The high
transitivity of the natural action of $S_n$ on $n$ objects enables
us to form the progenitors $\cal{P}$$_1$:=2$^{\star{n \choose
1}}\colon S_n$, $\cal{P}$$_2$:=2$^{\star{n \choose 2}}\colon S_n$
and $\cal{P}$$_3$:=2$^{\star{n \choose 3}}\colon S_n$.

Arguments similar to those used in the case $\cal{P}$$_1$ may be
applied in the other two cases naturally leading us to consider the
factored progenitors
\[ \frac{2^{\star {n \choose 2}}:S_n}{(t_{12}(23))^3}\mbox{ for }n\geq4\mbox{ and }\frac{2^{\star {n \choose 3}}:S_n}{(t_{123}(34))^3}\mbox{ for }n\geq6.\]

In all three cases the exceptional stabilizers and centralizers
encountered for small values of $n$ can be shown to lead straight to
interesting presentations of various finite groups \cite[Section 3.8]{PhD} but we
shall make no use of these results here.

\section{Coset Enumeration}\label{enumerate}

To prove that the homomorphic images under the relations appearing
in Theorem \ref{main} are finite we need to perform a double coset
enumeration to place an upper bound on the order of the target group
in each case.

The orders of all finite irreducible Coxeter groups, including those
of types $A_n$, $D_n$ and $E_n$, may be found listed in Humphreys
\cite[Table 2, p.44]{Coxetergroups}.

\subsection{$A_n$}
For $\cal{P}$$_1$ we enumerate the cosets by hand. Since
$t_it_j=(ij)t_i$ for $i,j\in\{1,\ldots,n\}$, $i\not=j$, any coset
representative must have length at most one. Since the stabilizer of a symmetric generator in our control group, $S_n^{(t_1)}$, clearly contains a subgroup isomorphic to
$S_{n-1}$ (namely the stabilizer in $S_n$ of the point 1). We have that $|S_n\colon S_n^{(t_1)}|\leq n$ and
$|S_n\colon S_n^{(\star)}|=1$, so the target group must contain the
image of $S_n$ to index at most $n+1$.

\subsection{$D_n$}\label{Dncosets}

 We shall prove:
\begin{Lem}\label{mainDn}
Let
\[ G:=\frac{2^{\star {n \choose 2}}:S_n}{(t_{12}(23))^3} \mbox{ for } n\geq4.\]
The representatives for the double cosets $S_nw S_n\subset G$ with
$w$ a word in the symmetric generators are [$\star$],
[$t_{12}$], [$t_{12}t_{34}$], $\ldots$, [$t_{12}t_{34}\ldots t_{2k-1,2k}$],
where $k$ is the largest integer such that $2k\leq n$. We thus have
$|G:S_n|\leq2^{n-1}$.
\end{Lem}
We shall prove this by using the following two lemmata.
\begin{Lem}\label{firstDn}
For the group $G$ as above, the double coset represented by the word
$t_{ab}\ldots t_{ij}\ldots t_{ik}\ldots t_{cd}$ may be
represented by a shorter word (ie if two symmetric generators in a given word have some index in common, then that word can be replaced by a shorter word).
\end{Lem}

\noindent{\bf Proof} The relation immediately tells us
$t_{12}t_{13}=(23)t_{12}$ and so $[t_{12}t_{13}]=[t_{12}]$, thus we can suppose our word has length at least three.
Using the high transitivity of the action of $S_n$ on $n$ points we may assume that
our word contains a subword of the form $t_{12}\ldots t_{34}t_{15}$ with no other occurrence of the index `1' and no other repetitions appearing anywhere between the symmetric generators $t_{12}$ and $t_{15}$ of this subword. Now,

\begin{eqnarray*}
t_{12}\ldots t_{34}t_{15}& = &t_{12}\ldots t_{34}t_{13}^2t_{15}\\
& = &t_{12}\ldots((14)t_{34})((35)t_{13})\\
& = &(14)(35)t_{24}\ldots t_{45}t_{13}\\
\end{eqnarray*}
and so the repeated indices can be `moved closer together'. Repeating the above, the two symmetric generators with the common index can eventually be placed side by side at which point our relation immediately shortens this word since $t_{12}t_{13}=(23)t_{12}$. Since our word has finite length we can easily repeat this procedure to eliminate all repetitions. \hspace{3.5mm}$\Box$

\begin{Lem}\label{secondDn}
$t_{12}t_{34}\sim t_{13}t_{24}$
\end{Lem}

\noindent{\bf Proof}
\[t_{12}t_{34}=t_{12}t_{34}t_{24}^2=t_{12}(23)t_{34}t_{24}=(23)t_{13}t_{34}t_{24}=(23)(14)t_{13}t_{24}\sim
t_{13}t_{24}.\]
$\Box$

\noindent{\bf Proof of Lemma \ref{mainDn}} By Lemma \ref{firstDn}
the indices appearing in any coset representative must be distinct.
By Lemma \ref{secondDn} the indices appearing in a word of length
two may be reordered. Since the indices are all distinct it follows
that the indices appearing in a coset representative of any length
may be reordered. The double cosets must therefore be [$\star$],
[$t_{12}$],$\ldots$,[$t_{12}\ldots t_{2k-1,2k}$],
where $k$ is the largest integer such that $2k\leq n$. There is therefore no more
than one double coset for each subset of $\{1,\ldots,n\}$ of even
size and so $|G: S_n|\leq 2^{n-1}$.\hspace{3.5mm}$\Box$

\subsection{E$_6$}

The coset enumeration in this case may also be performed by hand. We
list the cosets in Table I. Not every case is considered in this
table, however all remaining cases may be deduced from them as
follows. Since $t_{123}t_{145}\sim t_{124}t_{135}$ the S$_4$
permuting these indices ensures that for any three element subset
$\{a,b,c\}\subset\{1,\ldots,6\}$ the word
$t_{123}t_{145}t_{abc}$ will shorten. Since the only
non-collapsing word of length 3 is of the form
$t_{123}t_{456}t_{123}$ and $t_{123}t_{456}t_{123}\sim
t_{124}t_{356}t_{124}$ the S$_6$ permuting these indices
ensures that for any three element subset
$\{a,b,c\}\subset\{1,\ldots,6\}$ the word
$t_{123}t_{456}t_{123}t_{abc}$ will shorten and so all words
of length 4 shorten.

From this double coset enumeration we see that $|W($E$_6):$
S$_6|\leq1+20+30+20+1=72$. Our target group must therefore have
order at most $72\times|$S$_6|=51840$.
\begin{figure}[h]
\begin{center}
Table I: The coset enumeration for E$_6$
\begin{tabular}{|llc|}
\hline
Label [$w$]&Coset Stabilizing subgroup&$|N:N^{(w)}|$\\
\hline
[$\star$]&$N$&1\\
$\lbrack t_{123}\rbrack$&$N^{(t_{123})}\cong$ S$_3\times$S$_3$&20\\
$\lbrack t_{123}t_{145}\rbrack$&$N^{(t_{123}t_{145})}\cong$ S$_4$ since&30\\
&$t_{123}t_{145}=t_{123}t_{124}^2t_{145}\sim t_{123}(25)t_{124}$&\\
&$\sim t_{135}t_{124}$&\\
$\lbrack t_{123}t_{456}\rbrack$&$N^{(t_{123}t_{456})}\cong$ S$_3\times$S$_3$ since&20\\
$\lbrack t_{123}t_{456}t_{124}\rbrack$&$t_{123}t_{456}t_{124}=t_{123}t_{456}t_{145}^2t_{124}$&\\
=$\lbrack t_{356}t_{245}\rbrack$&\hspace{24mm}$=t_{123}(16)t_{456}(25)t_{145}$&\\
&$\hspace{24mm}\sim t_{356}t_{245}t_{145}$&\\
&$\hspace{24mm}\sim t_{356}t_{245}$&\\
$\lbrack t_{123}t_{456}t_{123}\rbrack$&$N^{(t_{123}t_{456}t_{123})}\cong$ S$_6$ since&1\\
&$t_{123}t_{456}t_{123}=t_{123}(34)t_{456}t_{356}t_{123}$&\\
&$\hspace{24mm}\sim t_{123}(34)t_{456}t_{356}t_{235}^2t_{123}$&\\
&$\hspace{24mm}=t_{124}t_{456}(26)t_{356}(15)t_{235}$&\\
&$\hspace{24mm}=t_{456}t_{124}t_{136}t_{235}$&\\
&$\hspace{24mm}=t_{456}t_{146}^2t_{124}t_{136}t_{235}$&\\
&$\hspace{24mm}=(15)t_{456}(62)t_{146}t_{136}t_{235}$&\\
&$\hspace{24mm}=t_{245}t_{146}t_{136}t_{235}$&\\
&$\hspace{24mm}=t_{245}(34)t_{146}t_{235}$&\\
&$\hspace{24mm}\sim t_{235}t_{146}t_{235}$&\\
\hline
\end{tabular}
\end{center}
\end{figure}

\subsection{E$_7$}

Since we expect both the index and the number of cosets to be much larger in this case than in the E$_6$ case (and in particular to be too unwieldy for a `by hand' approach to work) we use a computer, and in particular the algebra package {\sc Magma} \cite{MAGMA} to determine the index.

\begin{verbatim}
> S:=Sym(7);
> stab:=Stabilizer(S,{1,2,3});
> f,nn:=CosetAction(S,stab);
\end{verbatim}

Here we have defined a copy of the symmetric group S$_7$ (now named `\texttt{nn}') in its permutation representation defined by the action on the ${7 \choose3}=35$ subsets of cardinality 3 via the natural representation, and a homomorphism \texttt{f} from a copy of S$_7$ that acts on seven points to our new copy \texttt{nn}.

\begin{verbatim}
> 1^f(S!(3,4));
22
\end{verbatim}

The computer has labeled the set $\{1,2,3\}$ \texttt{1} and to find the label the computer has given to the set $\{1,2,4\}$ we find the image of \texttt{1} under the action of the permutation \texttt{f((1,2))}$\in$\texttt{nn}, finding that on this occasion the computer has given the set $\{1,2,4\}$ the label \texttt{22}.

\begin{verbatim}
> RR:=[<[1,22,1],f(S!(3,4))>];
> CT:=DCEnum(nn,RR,nn:Print:=5,Grain:=100);

Index: 576 = Rank: 10 = Edges: 40 = Status: Early closed = 
Time: 0.150
\end{verbatim}

The ordered sequence \texttt{RR} contains the sequence of symmetric generators $t_{123}t_{124}t_{123}$ and the permutation (34) that we are equating with this word to input our additional relation into the computer. The command \texttt{DCEnum} simply calls the double coset enumeration program of Bray and Curtis as described in \cite{cosetenumerator}.

The computer has found there to be at most 10 distinct double cosets
and that $|W($E$_7):$S$_7|\leq576$. Our target group must therefore
have order at most 576$\times|$S$_7|=2903040$.

\subsection{E$_8$}

Again, we use the computer to determine the index, each of the {\sc Magma} commands below being the same as those used in the previous section.

\begin{verbatim}
> S:=Sym(8);
> stab:=Stabilizer(S,{1,2,3});
> f,nn:=CosetAction(S,stab);
> 1^f(S!(3,4));
28
> RR:=[<[1,28,1],f(S!(3,4))>];
> CT:=DCEnum(nn,RR,nn:Print:=5,Grain:=100);
Index: 17280 = Rank: 35 = Edges: 256 = Status: Early closed = 
Time: 0.940
\end{verbatim}
\medskip
We see that $|W($E$_8):$S$_8|\leq17280$. Our target group must
therefore have order at most $17280\times|$S$_8|=696729600$.

\section{Representations}\label{reps}

In this section we use the symmetric presentations of Theorem 5 to
construct representations of the target groups and in doing so
verifying that we have the structures that we claim. In the $A_n$
and $D_n$ cases this is sufficient to show that the groups are what
we expect them to be.

\subsection{$W(A_n)$}

Since these groups are most naturally viewed as permutation groups
we shall construct the natural permutation representation. The
lowest degree of a permutation representation in which the control
group, $S_n$, acts faithfully is $n$, so the lowest degree of a
permutation representation in which the target group acts faithfully
is $n$. Since the control group already contains all possible
permutations of $n$ objects, the target group must be a permutation
group of at least $n+1$ objects. A permutation corresponding to a
symmetric generator must commute with its stabilizer in the
control group, namely $S_{n-1}$. There is only one such permutation
satisfying this: $t_i=(i,n+1)$. Since this has order two and
satisfies the relation we must therefore have that our target group
is isomorphic to $S_{n+1}\cong W(A_n)$.

\subsection{$W(D_n)$}

We shall use our symmetric generators to construct an elementary
Abelian 2-group lying outside our control group and thus verify that
our target group has structure $2^{n-1}:S_n$.

\begin{Lem}\label{somthing}
$t_{12}t_{34}=t_{34}t_{12}$
\end{Lem}

\noindent{\bf Proof}
\begin{eqnarray*}
t_{12}t_{34}t_{12}&=&t_{12}t_{34}t_{13}^2t_{12}\\
&=&t_{12}(14)t_{34}(23)t_{13}\\
&=&(14)(23)t_{34}t_{24}t_{13}\\
&=&(14)(t_{34}t_{24})t_{24}t_{13}\\
&=&t_{34}\\
\end{eqnarray*}\hspace{3.5mm}$\Box$

\begin{Lem}\label{elemab}
\label{6} The elements $e_{ij}:=(ij)t_{ij}$ for $1\leq i,j\leq n$ generate an elementary
Abelian 2-group.
\end{Lem}

\noindent{\bf Proof} Each of the element $e_{ij}$ have order 2 since the symmetric generators have order 2. If $i,j\notin\{k,l\}$ then by Lemma
\ref{somthing} $e_{ij}e_{kl}=e_{kl}e_{ij}$. Suppose $i=l$, then

\begin{eqnarray*}
e_{ij}e_{ik}e_{ij}e_{ik}& = &(ij)t_{ij}(ik)t_{ik}(ij)t_{ij}(ik)t_{ik}\\
& = &(ij)(ik)(ij)(ik)t_{ik}t_{ij}t_{jk}t_{ik}\\
& = &(ij)(ik)(ij)(ik)(jk)t_{ik}t_{jk}t_{ik}\\
& = &(ij)(ik)(ij)(ik)(jk)(ij)\\
& = &id_{S_n}
\end{eqnarray*}\hspace{3.5mm}$\Box$

\begin{Lem}\label{generators}
If $e_{ij}$ is as defined in Lemma \ref{elemab} then $e_{ij}e_{ik}=e_{jk}$ for $i\not=j\not=k\not=i$.
\end{Lem}

\noindent{\bf Proof} 
$e_{ij}e_{ik}=(ij)t_{ij}(ik)t_{ik}=(ij)(ik)t_{jk}t_{ik}=(ij)(ik)(ij)t_{jk}=(jk)t_{jk}=e_{jk}$ $\Box$

We have thus shown that there is an elementary Abelian
group of order $2^{n-1}$ lying outside the control group: the elements $e_{ij}$ defined in lemma \ref{elemab} each have order 2 (since the symmetric generators each have order 2), by lemma \ref{elemab} any two of the elements $e_{ij}$ commute and by lemma \ref{generators} the subgroup generated by these elements is clearly generated by the $n-1$ elements $e_{12},e_{13},\ldots, e_{1n}$. 

It is natural to represent the elements $e_{ij}$ as diagonal matrices
with -1 entries in the $i$ and $j$ positions. Using the natural
$n$-dimensional representation of $S_n$ as permutation matrices we
thus have been naturally lead to the following.

$$t_{12}=\left(\begin{tabular}{ccccc}
&-1&&&\\
-1&&&&\\
&&1&&\\
&&&$\ddots$&\\
&&&&1\\
\end{tabular}
\right)
$$

The
control group naturally acts on the group generated by the elements $e_{ij}$ by permuting the indices. In particular, recalling from the double coset enumeration of Section \ref{Dncosets} that $N$ has index at most $2^{n-1}$ in the target group, the above lemmas together show that our target group is isomorphic to the group $2^{n-1}:S_n\cong W(D_n)$.

\subsection{$W($E$_6$)}

In the case of E$_6$ we shall construct a 6 dimensional real
representation in which the control group acts as permutation
matrices. In such a representation the matrix corresponding to the
symmetric generator $t_{123}$ must:
\medskip
\begin{quote}
\begin{enumerate}
\item commute with the stabilizer of $t_{123}$;\\[-5pt]
\item have order two;\\[-5pt]
\item satisfy the relation.
\end{enumerate}
\end{quote}
\medskip
\noindent By condition 1 such a matrix be of the form
\medskip
$$
t_{123} = \left ( \begin{array}{c|c}  aI_3+bJ_3 \  & \ cJ_3  \cr
\hline c'J_3 \ & \ a'I_3+b'J_3
\end{array} \right )
$$

\medskip
\noindent where $I_3$ denotes the 3$\times$3 identity matrix and
$J_3$ denotes a 3$\times$3 matrix all the entries of which are 1.
Now, condition 2 requires
\[(aI_3+bJ_3)^2+3cc^\prime J_3=(a^\prime I_3+b^\prime
J_3)^2+3cc^\prime J=I_3\]implying that
\[c(a+a^\prime+3b+3b^\prime)=c^\prime(a+a^\prime+3b+3b^\prime)=0\]
\[a^2=a^{\prime2}=1\mbox{ and}\]
\[2ab+3b^2+3cc^\prime=2a^\prime b^\prime+3b^{\prime2}+3cc^\prime=0.\]
If our control group acts as permutation matrices then condition
3 implies that the determinant of the matrix for the symmetric
generators must be -1. This requires that
\[(a+3b)(a^\prime+3b^\prime)=-1.\] From these relations we are
naturally led to matrices of the form:

$$
t_{123} = \left ( \begin{array}{c|ccccc}  I_3-\frac{2}{3}J_3 & & &
\ \frac{1}{3}J_3 && \cr \hline 0_3  &&& \ I_3 &&
\end{array} \right )
$$

The representation of the control group we have used is not
irreducible and splits into two irreducible representations: the
subspace spanned by the vector $v$:=(1$^6$) and the subspace
$v^\perp$. The above matrices do not respect this decomposition
since they map $v$ to vectors of the form ($0^3$,1$^3$).
Consequently, the above representation of $W($E$_6)$ is irreducible.

\subsection{$W($E$_7$)}

Using arguments entirely analogous to those appearing in the
previous Section there is a 7 dimensional representation of $W($E$_7)$ in
which the control group acts as permutation matrices and we can
represent the symmetric generators for $W($E$_7)$ with matrices of
the form
$$
t_{123} = \left ( \begin{array}{c|ccccc}  I_3-\frac{2}{3}J_3 & & &
\ \frac{1}{3}J_{3\times4} && \cr \hline 0_{4\times3}  &&& \ I_4 &&
\end{array} \right )
$$
which again is irreducible.

\subsection{$W($E$_8$)}

Again using arguments entirely analogous to those used in the E$_6$
case there is an 8 dimensional representation of $W($E$_8)$ in which
the control group acts as permutation matrices and we can represent
the symmetric generators for E$_8$ with matrices of the form
$$
t_{1,2,3} = \left ( \begin{array}{c|ccccc}  I_3-\frac{2}{3}J_3 & & &
\ \frac{1}{3}J_{3\times5} && \cr \hline 0_{5\times3}  &&& \ I_5 &&
\end{array} \right )
$$
which again is irreducible.

\section{$\Bbb{Z}_2$ Representations of the groups $W(E_n)$}\label{modrep}

In this section we use the matrices obtained in Section \ref{reps}
for representing the Weyl groups of types E$_6$, E$_7$ and E$_8$ to exhibit representations
of these groups over $\Bbb{Z}_2$ and in doing so we identify the
structure of the groups in question.

\subsection{$W($E$_6$)}

Multiplying the matrices for our symmetric generators found in the
last Section by 3 ( $\equiv$1 (mod 2)) we find that these matrices,
working over $\Bbb{Z}_2$, are of the form:
$$
t_{123} = \left ( \begin{array}{c|c}  I_3 &   J_3  \cr \hline  0_3
&
 I_3
\end{array} \right )
$$

These matrices still satisfy the relation and the representation is
still irreducible for the same reason as in the real case as is
easily verified by {\sc Magma}. Consequently we see the isomorphism
$W(E_6$)$\cong$O$^-_6$(2):2 since all of our matrices preserve the
non-singular quadratic form $\sum_{i\not=j}x_ix_j$.

\subsection{$W($E$_7$)}

Similarly we obtain a representation of  $2\times$O$_7$(2) in the
E$_7$ case, accepting that the central involution must clearly act
trivially here. In this case the matrices preserve the non-singular
quadratic form defined by $xJ_7y^T$.

From the Atlas of Brauer Characters \cite[p.110]{ATLAS2} we see that
there is no irreducible $\Bbb{Z}_2$ representation of O$_7$(2) in 7
dimensions and this is precisely what we find here. The matrices for
the symmetric generators and the whole of the control group fix the
vector $v:=$(1$^7$). The space $v^\perp$ thus gives us a 6
dimensional $\Bbb{Z}_2$-module for this group to act on. It may be
easily verified with the aid of {\sc Magma} that this representation
is irreducible.

Since the above form is symplectic when restricted to this subspace
we immediately recover the classical exceptional isomorphism
O$_7$(2)$\cong$S$_6$(2).

(It is worth noting that in both the E$_6$ and E$_7$ cases the
symmetric generators may be interpreted as `bifid maps' acting on
the 27 lines of Schl\"{a}fli's general cubic surface and Hesse's 28
bitangents to the plane quartic curve respectively. See \cite[p.26
and p.46]{ATLAS} for details.)

\subsection{$W($E$_8$)}

Similarly we obtain a representation of 2\.{}O$_8^+$(2) in the E$_8$
case, again accepting that the central involution must clearly act
trivially. Like the E$_6$ case the matrices preserve the
non-singular quadratic form $\sum_{i\not=j}x_ix_j$.

Notice that working in an even number of dimensions removes the
irreducibility problem encountered with E$_7$ since the image of
(1$^8$) under the action of a symmetric generator is of the form
(0$^3$,1$^5$).\bigskip

{\bf Remark} Here we focused our attention on the simply laced
Coxeter groups. Analogous results may be obtained for other Coxeter
groups, but are much less enlightening. For example:
\[ \frac{2^{\star 2n}:W(\mbox{B}_{n-1})}{(t_1(12)(n+1,n+2))^3}\cong W(\mbox{B}_{n}) \]
\[\hspace{21mm}\frac{2^{\star n}:\mbox{S}_n}{(t_{1}(12))^5}\cong W(\mbox{H}_{n}) \mbox{ for }
n=3,4.\]\bigskip

{\bf Acknowledgments} I am deeply indebted to my PhD supervisor, Professor
Robert Curtis, for his continuing guidance and support throughout
this project without which this paper would not have been possible.
I am also grateful to Professor J\"{u}rgen M\"{u}ller for comments and
suggestions made about earlier versions of this paper.

\bibliographystyle{elsart-harv}

\end{document}